\documentclass[14pt,a4paper,fleqn]{article}
\usepackage{amsmath,amssymb,amsfonts}
\newtheorem{theorem}{Theorem}

\numberwithin{equation}{section}
\numberwithin{lemma}{section}
\numberwithin{theorem}{section}
\numberwithin{corollary}{section}

\allowdisplaybreaks

\begin{document}
\title{Finite summation formulas of  generalized  Kamp\'e de F\'eriet series}
\author{
Ashish Verma\footnote{Corresponding author}
\\ 
Department of Mathematics\\ Prof. Rajendra Singh (Rajju Bhaiya)\\ Institute of Physical Sciences for Study and Research \\  V. B. S Purvanchal University, Jaunpur  (U.P.)- 222003, India\\
vashish.lu@gmail.com}
\maketitle
\begin{abstract}
In the present  paper, author obtained finite summation  formulas for the  generalized  Kamp\'e de F\'eriet series. The peculiar outcome for  four  generalized  Lauricella  functions and confluent  forms of Lauricella series in $n$ variables are also drawn from the finite summation  formulas for the  generalized  Kamp\'e de F\'eriet series. \\[12pt]
\end{abstract}

\section{Introduction}

The multivariable generalization of  Kamp\'e de F\'eriet function  is defined as follows \cite{SK, SM} :
\begin{eqnarray}
&&F_{{l}:\,m_1;\dots;\,m_{n}}^{p:\,q_1;\dots;\,q_{n}}\Big[^{(a_p)\,\,:\,\,(b^{(1)}_{q_1})\,;\dots;\,(b^{(n)}_{q_{n}});}_{(\alpha_l) \,:\,(\beta^{(1)}_{m_1})\,;\dots;\,(\beta^{(n)}_{m_{n}});} \,\,x_1,\dots, x_n\Big]\notag\\ 
&&=\sum_{s_1, \dots, s_n= 0}^{\infty}\wedge(s_1,\,\dots,\,s_n)\prod_{i=1}^{n}\frac{{x_i}^{s_i}}{{s_i}{!}},\label{keq1}
\end{eqnarray}
where  
\begin{equation}
\wedge(s_1,\,\dots,\,s_n) = \frac{\displaystyle\prod_{j=1}^{p}(a_j)_{s_1+\cdots+s_n } \prod_{j=1}^{q_1}(b^{(1)}_j)_{s_1}\cdots\prod_{j=1}^{q_{n}}(b^{(n)}_{j})_{s_n}}{\displaystyle\prod_{j=1}^{l}(\alpha_j)_{s_1+\cdots+s_n } \prod_{j=1}^{m_1}(\beta^{(1)}_j)_{s_1}\cdots\prod_{j=1}^{m_{n}}(\beta^{(n)}_j)_{s_n}} \label{keq2},
\end{equation}
and, for convergence of (\ref{keq1}),
$1+l+m_{u}-p-q_{u}\geq0$, $ u=1,\dots,n$;
the equality holds when, in addition, either
$p>l $ and  when $|x_1|^{{1}/ {(p-l)}}+\dots+|x_n|^{{1}/ {(p-l)}}<1$,
or
$p\leq l$ and $\max({|x_1|,\dots,|x_n|})<1$.\\

Next,   recall the definition of derivative operator 
\begin{align*}
D_{x}f(x)=\lim_{h\to 0}\frac{f(x+h)-f(x)}{h},
\end{align*}
provided $f$ is differentiable at $x$. Also $D_{x}^{n}f{(x)}=D_{x}(D_{x}^{n-1}f(x))$,  $n=0,1,2,\dots.$\\

 Recently, Wang has established many infinite summation formulas of double hypergeometric functions \cite{XW}. Then Wang and Chen \cite{XW1} obtained  finite summation formulas of double hypergeometric functions by using certain summation theorems. Further, Sahai and Verma \cite{VS} derived finite summation formulas for the Srivastava's general triple hypergeometric function \cite{S} and the results of  Lauricella functions \cite{GL} and  Srivastava's triple hypergeometric functions \cite{S1, S2} are presented as example. These results
unified and generalized the  many results in \cite{XW1} for the three variable hypergeometric function. Motivated by their work, author present here several finite summation  formulas for the generalized  Kamp\'e de F\'eriet series. Certain particular cases leading to  finite summation formulas for four  generalized  Lauricella  functions and confluent  forms of Lauricella series in $n$ variables are also presented.\\
 
Following abbreviated notations are used. For example, write 
\begin{align}
&(a_{p}+k)=a_1+k,\dots,a_{p}+k,\notag\\
&(a^{i}_{p}+k)=a_1+k,\dots,a_{i-1}+k, a_{i+1}+k,\dots,a_{p}+k, \,\,\, i=1,\dots,p,\notag\\
&(b_{q_t}^{(t)}+k)=b_{1}^{(t)}+k, \dots\dots, b_{q_{t}}^{(t)}+k,\notag\\
&(b_{q_t}^{(t), i}+k)=b_{1}^{(t)}+k, \dots,b_{i-1}^{(t)}+k ,\, b_{i+1}^{(t)}+k,\dots, b_{q_{t}}^{(t)}+k, \,t=1,\dots,n,\,\, \,1\leq i \leq q_t.\end{align}
Also,  denote
\begin{align}
\qquad [a]_{k}=\prod_{j=1}^{p}\,(a_j)_{k},\qquad
[a^i]_{k}=\prod_{j=1, j\neq i}^{p}\,(a_j )_{k} ,\qquad [b^{(t)}]_{k}=\prod_{j=1}^{q_{t}}\,(b^{(t)}_j )_{k}, \notag\\
\qquad [b^{(t),i}]_{k}=\prod_{j=1, j\neq i}^{q_{t}}\,(b^{(t)}_j )_{k}\,\,,\qquad t=1,\dots, n,\,  \notag
\end{align}
where $k$  is non-negative integer and  $(a_j)_{k}$ is the Pochhammer symbol, \cite{R}.\\ 

\section{Finite summation formulas of generalized  Kamp\'e de F\'eriet series by derivative operator}
In this section,  author  obtained  the finite summation formulas of generalized  Kamp\'e de F\'eriet series by derivative operator. 
 The $r$th derivative on $x_1$ of generalized  Kamp\'e de F\'eriet series is obtained as follows:

\begin{align}
& D_{x_1}^{r}\{F_{{l}:\,m_1;\dots;\,m_{n}}^{p:\,q_1;\dots;\,q_{n}}\Big[^{(a_p)\,\,:\,\,(b^{(1)}_{q_1})\,;\dots;\,(b^{(n)}_{q_{n}});}_{(\alpha_l) \,:\,(\beta^{(1)}_{m_1})\,;\dots;\,(\beta^{(n)}_{m_{n}});} \,\,x_1,\dots, x_n\Big]\}\notag\\
&=\frac{[a]_{r}[b^{(1)}]_{r}}{[\alpha]_{r}[\beta^{(1)}]_{r}}\,\,F_{{l}:\,m_1;\dots;\,m_{n}}^{p:\,q_1;\dots;\,q_{n}}\Big[^{(a_p+r)\,\,:\,\,(b^{(1)}_{q_1}+r)\,;\,(b^{(2)}_{q_2});\dots;\,(b^{(n)}_{q_{n}});}_{(\alpha_l+r) \,:\,(\beta^{(1)}_{m_1}+r)\,;\,(\beta^{(2)}_{m_2})\,;\dots;\,(\beta^{(n)}_{m_{n}});} \,\,x_1,\dots, x_n\Big].\notag\\\label{s1}
\end{align}
By using the generalized Leibnitz formula
\begin{align*}
D_{x_1}^{r} \ \big(f(x_1) g(x_1)\big)=\sum_{k=0}^{r}{r\choose k} \ D_{x_1}^{r-k}f(x_1) \ D_{x_1}^{k}g(x_1)
\end{align*}
 and  (\ref{s1}),   derive the following finite summation formulas of generalized  Kamp\'e de F\'eriet series.
\begin{theorem}\label{t1}
The following finite summation formulas of generalized  Kamp\'e de F\'eriet series hold true:
\begin{align}
&\sum_{k=0}^{r}{r\choose k}\frac{[a^{i}]_{k}[b^{(1)}]_{k}}{[\alpha]_{k}[\beta^{(1)}]_{k}}\,x_{1}^{k}\,F_{{l}:\,m_1;\dots;\,m_{n}}^{p:\,q_1;\dots;\,q_{n}}\Big[^{(a_p+k)\,\,:\,\,(b^{(1)}_{q_1}+k)\,;\,(b^{(2)}_{q_2});\dots;\,(b^{(n)}_{q_{n}});}_{(\alpha_l+k) \,:\,(\beta^{(1)}_{m_1}+k)\,;\,(\beta^{(2)}_{m_2})\,;\dots;\,(\beta^{(n)}_{m_{n}});} \,\,x_1,\dots, x_n\Big]\notag\\
&=F_{{l}:\,m_1+1;m_{2};\dots;\,m_{n}}^{p:\,q_1+1;q_{2};\dots;\,q_{n}}\Big[^{(a_p)\,\,:\,a_{i}+r,\, (b^{(1)}_{q_1})\,;\,(b^{(2)}_{q_2});\dots;\,(b^{(n)}_{q_{n}});}_{(\alpha_l) \,:\,\,\,\,\,\,\,a_{i},(\beta^{(1)}_{m_1})\,\,;\,(\beta^{(2)}_{m_2})\,\dots;\,(\beta^{(n)}_{m_{n}});} \,\,x_1,\dots, x_n\Big],\label{a1s2}
\end{align}
where $i=1,\dots,p$;
\begin{align}
&\sum_{k=0}^{r}{r\choose k}\frac{[a]_{r}[b^{(1),i}]_{r}}{[\alpha]_{r}[\beta^{(1)}]_{r}}\,x_{1}^{k}\,F_{{l}:\,m_1;\dots;\,m_{n}}^{p:\,q_1;\dots;\,q_{n}}\Big[^{(a_p+k)\,\,:\,\,(b^{(1)}_{q_1}+k)\,;\,(b^{(2)}_{q_2});\dots;\,(b^{(n)}_{q_{n}});}_{(\alpha_l+k) \,:\,(\beta^{(1)}_{m_1}+k)\,;\,(\beta^{(2)}_{m_2})\,;\dots;\,(\beta^{(n)}_{m_{n}});} \,\,x_1,\dots, x_n\Big]\notag\\
&=F_{{l}:\,m_1;\dots;\,m_{n}}^{p:\,q_1;\dots;\,q_{n}}\Big[^{(a_p)\,\,:\,b^{(1)}_{i}+r,\, (b^{(1),i}_{q_1})\,;\,(b^{(2)}_{q_{2}});\dots;\,(b^{(n)}_{q_{n}});}_{(\alpha_l) \,:\,\,\,\,\,\,\,\,\,\,\,\,(\beta^{(1)}_{m_1})\,\,\,\,\,\,\,\,\,\,\,;\,(\beta^{(2)}_{m_2});\dots;\,(\beta^{(n)}_{m_{n}});} \,\,x_1,\dots, x_n\Big],\label{n2s2}
\end{align}
where $i=1,\dots,q_{1}$.

\end{theorem}
{\bf Proof:} For the proof of identity (\ref{a1s2}), using the  definition of generalized  Kamp\'e de F\'eriet series and the generalized Leibnitz formula for differentiation of a product of two functions, gives 
\begin{align*}
&D_{x_1}^{r}\{x_{1}^{a_{i}+r-1}F_{{l}:\,m_1;\dots;\,m_{n}}^{p:\,q_1;\dots;\,q_{n}}\Big[^{(a_p)\,\,:\,\,(b^{(1)}_{q_1})\,;\dots;\,(b^{(n)}_{q_{n}});}_{(\alpha_l) \,:\,(\beta^{(1)}_{m_1})\,;\dots;\,(\beta^{(n)}_{m_{n}});} \,\,x_1,\dots, x_n\Big]\}\\
&= \sum_{k=0}^{r}{r\choose k}D_{x_1}^{r-k}\{x_{1}^{a_{i}+r-1}\}D_{x_1}^{k}\{F_{{l}:\,m_1;\dots;\,m_{n}}^{p:\,q_1;\dots;\,q_{n}}\Big[^{(a_p)\,\,:\,\,(b^{(1)}_{q_1})\,;\dots;\,(b^{(n)}_{q_{n}});}_{(\alpha_l) \,:\,(\beta^{(1)}_{m_1})\,;\dots;\,(\beta^{(n)}_{m_{n}});} \,\,x_1,\dots, x_n\Big]\}\\
&=\,(a_{i})_{r}\,x_{1}^{a_{i}-1}\sum_{k=0}^{r}{r\choose k}\frac{[a^{i}]_{k}[b^{(1)}]_{k}}{[\alpha]_{k}[\beta^{(1)}]_{k}}\,x_{1}^{k}\\
&\qquad
\times F_{{l}:\,m_1;\dots;\,m_{n}}^{p:\,q_1;\dots;\,q_{n}}\Big[^{(a_p+k)\,\,:\,\,(b^{(1)}_{q_1}+k)\,;\,(b^{(2)}_{q_2});\dots;\,(b^{(n)}_{q_{n}});}_{(\alpha_l+k) \,:\,(\beta^{(1)}_{m_1}+k)\,;\,(\beta^{(2)}_{m_2})\,;\dots;\,(\beta^{(n)}_{m_{n}});} \,\,x_1,\dots, x_n\Big],
\end{align*}
using  (\ref{s1}) and some simplification in the second equality. Again,  combine $x_{1}^{a_{i}+r-1}$ with the variable $x_{1}$ in the  generalized  Kamp\'e de F\'eriet series and put the derivative operator $r$-times on $x_{1}$ to get the following result:
\begin{align*}
&D_{x_{1}}^{r}\{x_{1}^{a_{i}+r-1}F_{{l}:\,m_1;\dots;\,m_{n}}^{p:\,q_1;\dots;\,q_{n}}\Big[^{(a_p)\,\,:\,\,(b^{(1)}_{q_1})\,;\dots;\,(b^{(n)}_{q_{n}});}_{(\alpha_l) \,:\,(\beta^{(1)}_{m_1})\,;\dots;\,(\beta^{(n)}_{m_{n}});} \,\,x_1,\dots, x_n\Big]\}\\
&=\sum_{s_{1},\dots, s_{n}=0}^{\infty}\wedge(s_{1},\dots, s_{n}){(a_{i}+s_{1})_{r}\, x_{1}^{a_{i}-1}}\prod_{i=1}^{n}\frac{x_{i}^{s_{i}}}{s_{i}!}\\
&=(a_{i})_{r}\,x_{1}^{a_{i}-1}\,F_{{l}:\,m_1+1;m_{2};\dots;\,m_{n}}^{p:\,q_1+1;q_{2};\dots;\,q_{n}}\Big[^{(a_p)\,\,:\,a_{i}+r,\, (b^{(1)}_{q_1})\,;\,(b^{(2)}_{q_2});\dots;\,(b^{(n)}_{q_{n}});}_{(\alpha_l) \,:\,\,\,\,\,\,\,a_{i},(\beta^{(1)}_{m_1})\,\,;\,(\beta^{(2)}_{m_2})\,\dots;\,(\beta^{(n)}_{m_{n}});} \,\,x_1,\dots, x_n\Big].
\end{align*}
Equating the above two relations leads to (\ref{a1s2}).
 
The second result (\ref{n2s2}) are proved in a similar manner.

\begin{theorem}
The following finite summation formula of generalized  Kamp\'e de F\'eriet series holds true:
\begin{align}
&\sum_{k=0}^{r}{r\choose k}\frac{[a]_{k}[b^{(1)}]_{k}}{(\beta^{(1)}_{i}-r)_{k}[\alpha]_{k}[\beta^{(1)}]_{k}}\,x_{1}^{k}\,F_{{l}:\,m_1;\dots;\,m_{n}}^{p:\,q_1;\dots;\,q_{n}}\Big[^{(a_p+k)\,\,:\,\,(b^{(1)}_{q_1}+k)\,;\,(b^{(2)}_{q_2});\dots;\,(b^{(n)}_{q_{n}});}_{(\alpha_l+k) \,:\,(\beta^{(1)}_{m_1}+k)\,;\,(\beta^{(2)}_{m_2})\,;\dots;\,(\beta^{(n)}_{m_{n}});} \,\,x_1,\dots, x_n\Big]\notag\\
&=F_{{l}:\,m_1;\dots;\,m_{n}}^{p:\,q_1;\dots;\,q_{n}}\Big[^{(a_p)\,:\, \qquad(b^{(1)}_{q_1})\,\,\,\,\,\,\,\,\,\,;\,(b^{(2)}_{q_2});\dots;\,(b^{(n)}_{q_{n}});}_{(\alpha_l) \,:\,\beta^{(1)}_{i}-r,\,(\beta^{(1),i}_{m_1})\,;\,(\beta^{(2)}_{m_2})\,\dots;\,(\beta^{(n)}_{m_{n}});} \,\,x_1,\dots, x_n\Big],\label{1s2}
\end{align}
where $i=1,\dots,m_{1}$.
\end{theorem}
{\bf Proof:} Applying  the derivative operator on $x_{1}^{\beta^{(1)}_{i}-1}$$F_{{l}:\,m_1;\dots;\,m_{n}}^{p:\,q_1;\dots;\,q_{n}}(x_1,\dots,x_{n})$,  $r$-times, gives the formula in this theorem as explained in the proof of Theorem \ref{t1}. Omit the details.
\begin{theorem}
The following finite summation formulas of generalized  Kamp\'e de F\'eriet series hold true:
\begin{align}
&\sum_{k=0}^{r}{r\choose k}\frac{(-1)^{k}(1-\beta^{(1)}_{i})_{k}}{(2-\beta^{(1)}_{i}-r)_{k}}\,x_{1}^{k}\,F_{{l}:\,m_1;\dots;\,m_{n}}^{p:\,q_1;\dots;\,q_{n}}\Big[^{(a_p)\,:\, \qquad(b^{(1)}_{q_1})\,\,\,\,\,\,\,\,\,\,;\,(b^{(2)}_{q_2});\dots;\,(b^{(n)}_{q_{n}});}_{(\alpha_l) \,:\,\beta^{(1)}_{i}-k,\,(\beta^{(1),i}_{m_1})\,;\,(\beta^{(2)}_{m_2})\,\dots;\,(\beta^{(n)}_{m_{n}});} \,\,x_1,\dots, x_n\Big]\notag\\
&=\frac{(-1)^{r}[a]_{r}[b^{(1)}]_{r}}{(\beta^{(1)}_{i}-1)_{r}[\alpha]_{r}[\beta^{(1)}]_{r}}F_{{l}:\,m_1;\dots;\,m_{n}}^{p:\,q_1;\dots;\,q_{n}}\Big[^{(a_p+r)\,\,:\,\,(b^{(1)}_{q_1}+r)\,;\,(b^{(2)}_{q_2});\dots;\,(b^{(n)}_{q_{n}});}_{(\alpha_l+r) \,:\,(\beta^{(1)}_{m_1}+r)\,;\,(\beta^{(2)}_{m_2})\,;\dots;\,(\beta^{(n)}_{m_{n}});} \,\,x_1,\dots, x_n\Big],\label{s4}
\end{align}
\begin{align}
&\sum_{k=0}^{r}{r\choose k}\frac{(-1)^{k}(\beta^{(1)}_{i}+r-1)_{k}}{(\beta^{(1)}_{i})_{k}}\,x_{1}^{k}\,F_{{l}:\,m_1;\dots;\,m_{n}}^{p:\,q_1;\dots;\,q_{n}}\Big[^{(a_p)\,:\, \qquad(b^{(1)}_{q_1})\,\,\,\,\,\,\,\,\,\,;\,(b^{(2)}_{q_2});\dots;\,(b^{(n)}_{q_{n}});}_{(\alpha_l) \,:\,\beta^{(1)}_{i}+k,\,(\beta^{(1),i}_{m_1})\,;\,(\beta^{(2)}_{m_2})\,\dots;\,(\beta^{(n)}_{m_{n}});} \,\,x_1,\dots, x_n\Big]\notag\\
&=\frac{[a]_{r}[b^{(1)}]_{r}}{(\beta^{(1)}_{i}+r)_{r}[\alpha]_{r}[\beta^{(1)}]_{r}}x_{1}^{r}F_{{l}:\,m_1;\dots;\,m_{n}}^{p:\,q_1;\dots;\,q_{n}}\Big[^{(a_p+r)\,:\, \qquad(b^{(1)}_{q_1}+r)\,\,\,\,\,\,\,\,\,\,;\,(b^{(2)}_{q_2});\dots;\,(b^{(n)}_{q_{n}});}_{(\alpha_l+r) \,:\,\beta^{(1)}_{i}+2r,\,(\beta^{(1),i}_{m_1}+r)\,;\,(\beta^{(2)}_{m_2})\,\dots;\,(\beta^{(n)}_{m_{n}});} \,\,x_1,\dots, x_n\Big],\label{s5}
\end{align}
where $i=1,\dots,m_{1}$.
\end{theorem}
{\bf Proof:} First prove identity (\ref{s4}). 
 From the definition of generalized  Kamp\'e de F\'eriet series and  the generalized Leibnitz formula for differentiation of a product of two functions, gives the following result:
\begin{align*}
&D_{x_{1}}^{r}\{x_{1}^{1-\beta^{(1)}_{i}}\times x_{1}^{1-\beta^{(1)}_{i}}F_{{l}:\,m_1;\dots;\,m_{n}}^{p:\,q_1;\dots;\,q_{n}}\Big[^{(a_p)\,\,:\,\,(b^{(1)}_{q_1})\,;\dots;\,(b^{(n)}_{q_{n}});}_{(\alpha_l) \,:\,(\beta^{(1)}_{m_1})\,;\dots;\,(\beta^{(n)}_{m_{n}});} \,\,x_1,\dots, x_n\Big]\}\\
&=\sum_{k=0}^{r}{r\choose k}D_{x_{1}}^{r-k}\{x_{1}^{1-\beta^{(1)}_{i}}\}D_{x_1}^{k}\{x_{1}^{\beta^{(1)}_{i}-1}\,F_{{l}:\,m_1;\dots;\,m_{n}}^{p:\,q_1;\dots;\,q_{n}}\Big[^{(a_p)\,\,:\,\,(b^{(1)}_{q_1})\,;\dots;\,(b^{(n)}_{q_{n}});}_{(\alpha_l) \,:\,(\beta^{(1)}_{m_1})\,;\dots;\,(\beta^{(n)}_{m_{n}});} \,\,x_1,\dots, x_n\Big]\}\\
&=\sum_{k=0}^{r}(-1)^{r+k}{r\choose k}\frac{(\beta^{(1)}_{i}-1)_{r}(1-\beta^{(1)}_{i})_{k}}{(2-\beta^{(1)}_{i}-r)_{k}\,x_{1}^{r}}\\
&\quad\times F_{{l}:\,m_1;\dots;\,m_{n}}^{p:\,q_1;\dots;\,q_{n}}\Big[^{(a_p)\,:\, \qquad(b^{(1)}_{q_1})\,\,\,\,\,\,\,\,\,\,;\,(b^{(2)}_{q_2});\dots;\,(b^{(n)}_{q_{n}});}_{(\alpha_l) \,:\,\beta^{(1)}_{i}-k,\,(\beta^{(1),i}_{m_1})\,;\,(\beta^{(2)}_{m_2})\,\dots;\,(\beta^{(n)}_{m_{n}});} \,\,x_1,\dots, x_n\Big].
\end{align*}
Now using the derivative operator on generalized  Kamp\'e de F\'eriet series for $r$-times directly and equating with the above equality gives (\ref{s4}) after some simplification. Next, applying the operator $D_{x_{1}}^{r}$ on 
\begin{align*}
x_{1}^{1-\beta^{(1)}_{i}-r}\times x_{1}^{\beta^{(1)}_{i}+r-1}F_{{l}:\,m_1;\dots;\,m_{n}}^{p:\,q_1;\dots;\,q_{n}}\Big[^{(a_p)\,:\, \qquad(b^{(1)}_{q_1})\,\,\,\,\,\,\,\,\,\,;\,(b^{(2)}_{q_2});\dots;\,(b^{(n)}_{q_{n}});}_{(\alpha_l) \,:\,\beta^{(1)}_{i}+r,\,(\beta^{(1),i}_{m_1})\,;\,(\beta^{(2)}_{m_2})\,\dots;\,(\beta^{(n)}_{m_{n}});} \,\,x_1,\dots, x_n\Big],
\end{align*}
and proceeding as in the proof of (\ref{s4}) gives identity (\ref{s5}).
\begin{theorem}\label{t2}
The following finite summation formulas of generalized  Kamp\'e de F\'eriet series hold true:
\begin{align}
&\sum_{k=0}^{r}\frac{(-r)_{k}}{(a_{i}-r+1)_{k}}\,\,F_{{l+1}:\,m_1;\dots;\,m_{n}}^{p+1:\,q_1;\dots;\,q_{n}}\Big[^{1+k, (a_p)\,\,:\,\,(b^{(1)}_{q_1})\,;\dots;\,(b^{(n)}_{q_{n}});}_{1, (\alpha_l) \,:\,(\beta^{(1)}_{m_1})\,;\dots;\,(\beta^{(n)}_{m_{n}});} \,\,\frac{1}{x_1},\dots, \frac{1}{x_1}\Big]\notag\\
&=\,\frac{a_{i}-r}{a_{i}}F_{{l+1}:\,m_1;\dots;\,m_{n}}^{p+1:\,q_1;\dots;\,q_{n}}\Big[^{1-a_i+r, (a_p)\,\,:\,\,(b^{(1)}_{q_1})\,;\dots;\,(b^{(n)}_{q_{n}});}_{1-a_i, (\alpha_l) \,:\,(\beta^{(1)}_{m_1})\,;\dots;\,(\beta^{(n)}_{m_{n}});} \,\,\frac{1}{x_1},\dots, \frac{1}{x_1}\Big],\label{1s6}
\end{align}
where $i=1,\dots, p;$
\begin{align}
&\sum_{k=0}^{r}\frac{(-r)_{k}}{(b_{i}^{(1)}-r+1)_{k}}\,\,F_{{l}:\,m_1+1; m_2; \dots;\,m_{n}}^{p:\,q_1+1; q_2; \dots;\,q_{n}}\Big[^{ (a_p)\,\,:1+k, (b^{(1)}_{q_1})\,;\dots;\,(b^{(n)}_{q_{n}});}_{ (\alpha_l) \,:1, (\beta^{(1)}_{m_1})\,;\dots;\,(\beta^{(n)}_{m_{n}});} \,\,\frac{1}{x_1}, x_{2, }\dots, {x_n}\Big]\notag\\
&=\,\frac{b_{i}^{(1)}-r}{b_{i}^{(1)}}F_{{l}:\,m_1+1; m_2; \dots;\,m_{n}}^{p:\,q_1+1; q_2;\dots;\,q_{n}}\Big[^{ (a_p)\,\,:1-b_{i}^{(1)}+r , (b^{(1)}_{q_1})\,;\dots;\,(b^{(n)}_{q_{n}});}_{ (\alpha_l) \,:1-b_{i}^{(1)}, (\beta^{(1)}_{m_1})\,;\dots;\,(\beta^{(n)}_{m_{n}});} \,\,\frac{1}{x_1}, x_{2}, \dots, {x_n}\Big],\label{m2s6}
\end{align}
where $i=1,\dots, q_1;$

\end{theorem}
{\bf Proof:} To proof identity (\ref{1s6}). Identity (\ref{m2s6})  can be proved in an analogous manner. Calculating $r$th derivatives on $x_1$ of $x_{1}^{a_{i}-1} F_{{l}:\,m_1;\dots;\,m_{n}}^{p:\,q_1;\dots;\,q_{n}}(\frac{1}{x_1},\dots, \frac{1}{x_{n}})$,  gives
\begin{align*}
&D_{x_{1}}^{r}\{x_{1}^{a_{i}-1}F_{{l}:\,m_1;\dots;\,m_{n}}^{p:\,q_1;\dots;\,q_{n}}\Big[^{(a_p)\,\,:\,\,(b^{(1)}_{q_1})\,;\dots;\,(b^{(n)}_{q_{n}});}_{(\alpha_l) \,:\,(\beta^{(1)}_{m_1})\,;\dots;\,(\beta^{(n)}_{m_{n}});} \,\frac{1}{x_1},\dots, \frac{1}{x_n}\Big]\}\\
&=(-1)^{r} \frac{(1-a_{i})_{r}}{x_{1}^{r+1-a_{i}}}\,F_{{l+1}:\,m_1;\dots;\,m_{n}}^{p+1:\,q_1;\dots;\,q_{n}}\Big[^{1-a_i+r, (a_p)\,\,:\,\,(b^{(1)}_{q_1})\,;\dots;\,(b^{(n)}_{q_{n}});}_{1-a_i, (\alpha_l) \,:\,(\beta^{(1)}_{m_1})\,;\dots;\,(\beta^{(n)}_{m_{n}});} \,\,\frac{1}{x_1},\dots, \frac{1}{x_1}\Big].
\end{align*}
Alternatively, by the generalized Leibnitz formula for differentiation of a product of two functions,  get the following result:
\begin{align*}
&D_{x_1}^{r}\{x_{1}^{a_{i}-1}F_{{l}:\,m_1;\dots;\,m_{n}}^{p:\,q_1;\dots;\,q_{n}}\Big[^{(a_p)\,\,:\,\,(b^{(1)}_{q_1})\,;\dots;\,(b^{(n)}_{q_{n}});}_{(\alpha_l) \,:\,(\beta^{(1)}_{m_1})\,;\dots;\,(\beta^{(n)}_{m_{n}});} \,\frac{1}{x_1},\dots, \frac{1}{x_n}\Big]\}\\
&=\sum_{k=0}^{r}{r\choose k}\, D_{x_1}^{r-k}\{x_{1}^{a_{i}}\}\,D_{x_1}^{k}\{x_{1}^{-1}F_{{l}:\,m_1;\dots;\,m_{n}}^{p:\,q_1;\dots;\,q_{n}}\Big[^{(a_p)\,\,:\,\,(b^{(1)}_{q_1})\,;\dots;\,(b^{(n)}_{q_{n}});}_{(\alpha_l) \,:\,(\beta^{(1)}_{m_1})\,;\dots;\,(\beta^{(n)}_{m_{n}});} \,\frac{1}{x_1},\dots, \frac{1}{x_n}\Big]\}\\
&=(-1)^{r} (-a_{i})_{r}\, x_{1}^{a_{i}-r-1}\sum_{k=0}^{r}\frac{(-r)_{k}}{(a_{i}-r+1)_{k}}F_{{l+1}:\,m_1;\dots;\,m_{n}}^{p+1:\,q_1;\dots;\,q_{n}}\Big[^{1+k, (a_p)\,\,:\,\,(b^{(1)}_{q_1})\,;\dots;\,(b^{(n)}_{q_{n}});}_{1, (\alpha_l) \,:\,(\beta^{(1)}_{m_1})\,;\dots;\,(\beta^{(n)}_{m_{n}});} \,\,\frac{1}{x_1},\dots, \frac{1}{x_1}\Big].
\end{align*}
Equating the above two identities and after some simplification,  get (\ref{1s6}). 
This completes the proof of the theorem.
\begin{theorem}The following finite summation formulas of generalized  Kamp\'e de F\'eriet series hold true:
\begin{align}
&\sum_{k=0}^{r}\frac{(-r)_{k}}{(2-a_{i}-r)_{k}}\,\,F_{{l+1}:\,m_1;\dots;\,m_{n}}^{p+1:\,q_1;\dots;\,q_{n}}\Big[^{1+k, (a_p)\,\,:\,\,(b^{(1)}_{q_1})\,;\dots;\,(b^{(n)}_{q_{n}});}_{1, (\alpha_l) \,:\,(\beta^{(1)}_{m_1})\,;\dots;\,(\beta^{(n)}_{m_{n}});} \,\,\frac{1}{x_1},\dots, \frac{1}{x_1}\Big]\notag\\
&=\,\frac{a_{i}+r-1}{a_{i}-1}F_{{l}:\,m_1;\dots;\,m_{n}}^{p:\,q_1;\dots;\,q_{n}}\Big[^{a_i+r, (a_p^{i})\,\,:\,\,(b^{(1)}_{q_1})\,;\dots;\,(b^{(n)}_{q_{n}});}_{ (\alpha_l) \,:\,(\beta^{(1)}_{m_1})\,;\dots;\,(\beta^{(n)}_{m_{n}});} \,\,\frac{1}{x_1},\dots, \frac{1}{x_1}\Big],\label{1s7}
\end{align}
where $i=1,\dots, p;$
\begin{align}
&\sum_{k=0}^{r}\frac{(-r)_{k}}{(2-b_{i}^{(1)}-r)_{k}}\,\,F_{{l+1}:\,m_1;\dots;\,m_{n}}^{p+1:\,q_1;\dots;\,q_{n}}\Big[^{ (a_p)\,\,:1+k, (b^{(1)}_{q_1})\,;\dots;\,(b^{(n)}_{q_{n}});}_{ (\alpha_l) \,:1, (\beta^{(1)}_{m_1})\,;\dots;\,(\beta^{(n)}_{m_{n}});} \,\,\frac{1}{x_1}, x_{2, }\dots, {x_n}\Big]\notag\\
&=\,\frac{b_{i}^{(1)}+r-1}{b_{i}^{(1)}-1}F_{{l}:\,m_1;\dots;\,m_{n}}^{p:\,q_1;\dots;\,q_{n}}\Big[^{ (a_p)\,\,: b_{i}^{(1)}+r , (b^{(1), i}_{q_1})\,;\dots;\,(b^{(n)}_{q_{n}});}_{ (\alpha_l) \,: (\beta^{(1)}_{m_1})\,;\dots;\,(\beta^{(n)}_{m_{n}});} \,\,\frac{1}{x_1}, x_{2}, \dots, {x_n}\Big],\label{m2s8}
\end{align}
where $i=1,\dots, q_1;$

\end{theorem}
{\bf Proof:}   Outline the proof of (\ref{1s7}). Multiplying $x_{1}^{1-a_{i}}\times x_{1}^{-1}$  on the left-hand side of $F_{{l}:\,m_1;\dots;\,m_{n}}^{p:\,q_1;\dots;\,q_{n}}(\frac{1}{x_1},\dots, \frac{1}{x_{n}})$ series, and then using the derivative operator on them as in the proof of (\ref{1s6}) gives (\ref{1s7}). Identity (\ref{m2s8})  are proved in a similar manner.
\begin{theorem}The following finite summation formulas of generalized  Kamp\'e de F\'eriet series hold true:
\begin{align}
&\sum_{k=0}^{r}\frac{(-r)_{k}}{(2-2r)_{k}}\,\,F_{{l+1}:\,m_1;\dots;\,m_{n}}^{p+1:\,q_1;\dots;\,q_{n}}\Big[^{1+k, (a_p)\,\,:\,\,(b^{(1)}_{q_1})\,;\dots;\,(b^{(n)}_{q_{n}});}_{1, (\alpha_l) \,:\,(\beta^{(1)}_{m_1})\,;\dots;\,(\beta^{(n)}_{m_{n}});} \,\,\frac{1}{x_1},\dots, \frac{1}{x_1}\Big]\notag\\
&=\,\frac{2r-1}{r-1}F_{{l+1}:\,m_1;\dots;\,m_{n}}^{p+1:\,q_1;\dots;\,q_{n}}\Big[^{2r, (a_p)\,\,:\,\,(b^{(1)}_{q_1})\,;\dots;\,(b^{(n)}_{q_{n}});}_{r, (\alpha_l) \,:\,(\beta^{(1)}_{m_1})\,;\dots;\,(\beta^{(n)}_{m_{n}});} \,\,\frac{1}{x_1},\dots, \frac{1}{x_1}\Big],\label{1s9}
\end{align}
where $i=1,\dots, p;$
\begin{align}
&\sum_{k=0}^{r}\frac{(-r)_{k}}{(2-2r)_{k}}\,\,F_{{l}:\,m_1+1; m_2;\dots;\,m_{n}}^{p:\,q_1+1; q_2;\dots;\,q_{n}}\Big[^{ (a_p)\,\,:1+k, (b^{(1)}_{q_1})\,;\dots;\,(b^{(n)}_{q_{n}});}_{ (\alpha_l) \,:1, (\beta^{(1)}_{m_1})\,;\dots;\,(\beta^{(n)}_{m_{n}});} \,\,\frac{1}{x_1}, x_{2, }\dots, {x_n}\Big]\notag\\
&=\,\frac{2r-1}{r-1}F_{{l}:\,m_1+1; m_2;\dots;\,m_{n}}^{p:\,q_1+1; q_2;\dots;\,q_{n}}\Big[^{ (a_p)\,\,: 2r , (b^{(1)}_{q_1})\,;\dots;\,(b^{(n)}_{q_{n}});}_{ (\alpha_l) \,: r, (\beta^{(1)}_{m_1})\,;\dots;\,(\beta^{(n)}_{m_{n}});} \,\,\frac{1}{x_1}, x_{2}, \dots, {x_n}\Big],\label{m2s10}
\end{align}
where $i=1,\dots, q_1.$

\end{theorem}
{\bf Proof:} Here give an outline of proof of (\ref{1s9}). Multiplying $x_{1}^{1-r}\times x_{1}^{-1}$  on the left-hand side of $F_{{l}:\,m_1;\dots;\,m_{n}}^{p:\,q_1;\dots;\,q_{n}}(\frac{1}{x_1},\dots, \frac{1}{x_{n}})$-series, and then applying the derivative operator on them as in the proof of (\ref{1s6}) gives (\ref{1s9}). Second result (\ref{m2s10})  are proved in an analogous manner.
\begin{theorem}\label{t3}
The following finite summation formulas of generalized  Kamp\'e de F\'eriet series hold true:
\begin{align}
&\sum_{k=0}^{r}\frac{(-r)_{k}}{(1-2r)_{k}}\,\,F_{{l+1}:\,m_1;\dots;\,m_{n}}^{p+1:\,q_1;\dots;\,q_{n}}\Big[^{1+k, (a_p)\,\,:\,\,(b^{(1)}_{q_1})\,;\dots;\,(b^{(n)}_{q_{n}});}_{1, (\alpha_l) \,:\,(\beta^{(1)}_{m_1})\,;\dots;\,(\beta^{(n)}_{m_{n}});} \,\,\frac{1}{x_1},\dots, \frac{1}{x_1}\Big]\notag\\
&=\,2\,F_{{l+1}:\,m_1;\dots;\,m_{n}}^{p+1:\,q_1;\dots;\,q_{n}}\Big[^{1+2r, (a_p)\,\,:\,\,(b^{(1)}_{q_1})\,;\dots;\,(b^{(n)}_{q_{n}});}_{1+r, (\alpha_l) \,:\,(\beta^{(1)}_{m_1})\,;\dots;\,(\beta^{(n)}_{m_{n}});} \,\,\frac{1}{x_1},\dots, \frac{1}{x_1}\Big],\label{1s11}
\end{align}
\begin{align}
&\sum_{k=0}^{r}\frac{(-r)_{k}}{(1-2r)_{k}}\,\,F_{{l}:\,m_1+1; m_2;\dots;\,m_{n}}^{p:\,q_1+1; q_2;\dots;\,q_{n}}\Big[^{ (a_p)\,\,:1+k, (b^{(1)}_{q_1})\,;\dots;\,(b^{(n)}_{q_{n}});}_{ (\alpha_l) \,:1, (\beta^{(1)}_{m_1})\,;\dots;\,(\beta^{(n)}_{m_{n}});} \,\,\frac{1}{x_1}, x_{2, }\dots, {x_n}\Big]\notag\\
&=\,2\,F_{{l}:\,m_1+1; m_2;\dots;\,m_{n}}^{p:\,q_1+1; q_2;\dots;\,q_{n}}\Big[^{ (a_p)\,\,: 1+2r , (b^{(1)}_{q_1})\,;\dots;\,(b^{(n)}_{q_{n}});}_{ (\alpha_l) \,: 1+r, (\beta^{(1)}_{m_1})\,;\dots;\,(\beta^{(n)}_{m_{n}});} \,\,\frac{1}{x_1}, x_{2}, \dots, {x_n}\Big].\label{m2s12}
\end{align}
\end{theorem}
{\bf Proof:}  Applying the derivative operator  $r$-times as follows, gives 
\begin{align*}
&D_{x_1}^{r}\{x_{1}^{-r-1}F_{{l}:\,m_1;\dots;\,m_{n}}^{p:\,q_1;\dots;\,q_{n}}\Big[^{(a_p)\,\,:\,\,(b^{(1)}_{q_1})\,;\dots;\,(b^{(n)}_{q_{n}});}_{(\alpha_l) \,:\,(\beta^{(1)}_{m_1})\,;\dots;\,(\beta^{(n)}_{m_{n}});} \,\frac{1}{x_1},\dots, \frac{1}{x_1}\Big]\}\\
&=\sum_{s_1,\dots, s_n=0}^{\infty}\wedge(s_1,\dots, s_n)\, \frac{x_{1}^{-s_1-s_2\cdots-s_n-2r-1}}{\prod_{i=1}^{n}s_{i}!}
(-1)^{r}\,(s_1+\cdots+s_n+r+1)_{r}\\
&= (-1)^{r}\,(1+r)_{r}\, x_{1}^{-2r-1}F_{{l+1}:\,m_1;\dots;\,m_{n}}^{p+1:\,q_1;\dots;\,q_{n}}\Big[^{1+2r, (a_p)\,\,:\,\,(b^{(1)}_{q_1})\,;\dots;\,(b^{(n)}_{q_{n}});}_{1+r, (\alpha_l) \,:\,(\beta^{(1)}_{m_1})\,;\dots;\,(\beta^{(n)}_{m_{n}});} \,\,\frac{1}{x_1},\dots, \frac{1}{x_1}\Big].
\end{align*}
Again, put the generalized Leibnitz formula for differentiation of the product of the following two functions and establish
\begin{align*}
&D_{x_1}^{r}\{x_{1}^{-r-1}F_{{l}:\,m_1;\dots;\,m_{n}}^{p:\,q_1;\dots;\,q_{n}}\Big[^{(a_p)\,\,:\,\,(b^{(1)}_{q_1})\,;\dots;\,(b^{(n)}_{q_{n}});}_{(\alpha_l) \,:\,(\beta^{(1)}_{m_1})\,;\dots;\,(\beta^{(n)}_{m_{n}});} \,\frac{1}{x_1},\dots, \frac{1}{x_1}\Big]\}\\
&=\sum_{k=0}^{r}{r\choose k}\,D_{x_1}^{r-k}\{x_1^{-r}\}D_{x_1}^{k}\{x_1^{-1}F_{{l}:\,m_1;\dots;\,m_{n}}^{p:\,q_1;\dots;\,q_{n}}\Big[^{(a_p)\,\,:\,\,(b^{(1)}_{q_1})\,;\dots;\,(b^{(n)}_{q_{n}});}_{(\alpha_l) \,:\,(\beta^{(1)}_{m_1})\,;\dots;\,(\beta^{(n)}_{m_{n}});} \,\frac{1}{x_1},\dots, \frac{1}{x_1}\Big]\}\\
&=(-1)^{r} (r)_{r}\,x_{1}^{-2r-1}\sum_{k=0}^{r}\frac{(-r)_{k}}{(1-2r)_{k}}F_{{l+1}:\,m_1;\dots;\,m_{n}}^{p+1:\,q_1;\dots;\,q_{n}}\Big[^{1+k, (a_p)\,\,:\,\,(b^{(1)}_{q_1})\,;\dots;\,(b^{(n)}_{q_{n}});}_{1, (\alpha_l) \,:\,(\beta^{(1)}_{m_1})\,;\dots;\,(\beta^{(n)}_{m_{n}});} \,\,\frac{1}{x_1},\dots, \frac{1}{x_1}\Big].
\end{align*}
Equating the above two equalities,  gives the finite summation formula (\ref{1s11}). Result (\ref{m2s12}) are proved in an analogous manner. This completes the proof of this theorem.
\begin{theorem}
The following finite summation formulas of generalized  Kamp\'e de F\'eriet series hold true:
\begin{align}
&\sum_{k=0}^{r}{r\choose k}\frac{(-1)^{k}(-r)_{k}}{(1+a_{i}-r)_{k}}\,\,F_{{l+\alpha}:\,m_1;\dots;\,m_{n}}^{p+\alpha:\,q_1;\dots;\,q_{n}}\Big[^{\,\,\,\,\,\frac{1+r}{\alpha},\dots, \frac{\alpha+r}{\alpha},\,\, (a_p)\,\,:\,\,(b^{(1)}_{q_1})\,;\dots;\,(b^{(n)}_{q_{n}});}_{\frac{1+r-k}{\alpha},\dots, \frac{\alpha+r-k}{\alpha}, (\alpha_l) \,:\,(\beta^{(1)}_{m_1})\,;\dots;\,(\beta^{(n)}_{m_{n}});} \,x_1^{\alpha},\dots, x_1^{\alpha}\Big]\notag\\
&=\frac{(-1)^{r}(1+a_{i})_{r}}{(-a_{i})_{r}}\,F_{{l+\alpha}:\,m_1;\dots;\,m_{n}}^{p+\alpha:\,q_1;\dots;\,q_{n}}\Big[^{\frac{1+a_{i}+r}{\alpha},\dots, \frac{\alpha+a_{i}+r}{\alpha}, (a_p)\,\,:\,\,(b^{(1)}_{q_1})\,;\dots;\,(b^{(n)}_{q_{n}});}_{\,\frac{1+a_{i}}{\alpha},\dots, \frac{\alpha+a_{i}}{\alpha},\, (\alpha_l) \,:\,(\beta^{(1)}_{m_1})\,;\dots;\,(\beta^{(n)}_{m_{n}});} \,x_1^{\alpha},\dots, x_1^{\alpha}\Big],\,\,\,\alpha\geq 2,\label{1s13}
\end{align}
where $i=1,\dots, p;$
\begin{align}
&\sum_{k=0}^{r}{r\choose k}\frac{(-1)^{k}(-r)_{k}}{(1+b_{i}^{(1)}-r)_{k}}\,\,F_{{l}:\,m_1+\alpha; m_2;\dots;\,m_{n}}^{p:\,q_1+\alpha; q_2;\dots;\,q_{n}}\Big[^{ (a_p)\,\,:\,\,\,\,\,\frac{1+r}{\alpha},\dots, \frac{\alpha+r}{\alpha},\, (b^{(1)}_{q_1})\,;\dots;\,(b^{(n)}_{q_{n}});}_{ (\alpha_l) \,:\frac{1+r-k}{\alpha},\dots, \frac{\alpha+r-k}{\alpha}, (\beta^{(1)}_{m_1})\,;\dots;\,(\beta^{(n)}_{m_{n}});} \,\,{x_1^{\alpha}}, x_{2, }\dots, {x_n}\Big]\notag\\
&=\frac{(-1)^{r}(1+b_{i}^{(1)})_{r}}{(-b_{i}^{(1)})_{r}}\,F_{{l}:\,m_1+\alpha; m_2;\dots;\,m_{n}}^{p:\,q_1+\alpha; q_2;\dots;\,q_{n}}\Big[^{ (a_p)\,\,: \frac{1+b_{i}^{(1)}+r}{\alpha},\dots, \frac{\alpha+b_{i}^{(1)}+r}{\alpha}, (b^{(1)}_{q_1})\,;\dots;\,(b^{(n)}_{q_{n}});}_{ (\alpha_l) \,: \frac{1+b_{i}^{(1)}}{\alpha},\dots, \frac{\alpha+b_{i}^{(1)}}{\alpha},\, (\beta^{(1)}_{m_1})\,;\dots;\,(\beta^{(n)}_{m_{n}});} \,{x_1^{\alpha}}, x_{2}, \dots, {x_n}\Big],\,\,\,\alpha\geq 2,\label{m2s14}
\end{align}
where $i=1,\dots, q_1.$
\end{theorem}
{\bf Proof:} First prove (\ref{1s13}). Using the definition of derivative operator on the following product $r$-times, gives 
\begin{align*}
&D_{x_1}^{r}\{x_{1}^{a_{i}+r}F_{{l}:\,m_1;\dots;\,m_{n}}^{p:\,q_1;\dots;\,q_{n}}\Big[^{(a_p)\,\,:\,\,(b^{(1)}_{q_1})\,;\dots;\,(b^{(n)}_{q_{n}});}_{(\alpha_l) \,:\,(\beta^{(1)}_{m_1})\,;\dots;\,(\beta^{(n)}_{m_{n}});} \,x_1^{\alpha},\dots, x_1^{\alpha}\Big]\}\\
&=\sum_{s_1,\dots,s_n=0}^{\infty}\wedge(s_1,\dots,s_n)\frac{x_{1}^{\alpha (s_1+\cdots+s_n)+a_{i}}}{\prod_{i=1}^{n}s_{i}!}\frac{(\alpha(s_1+\cdots+s_n) +a_{i}+r)!}{(\alpha(s_1+\cdots+s_n)+a_{i})!}\\
&=x_{1}^{a_{i}}\sum_{s_1,\dots,s_n=0}^{\infty}\wedge(s_1,\dots,s_n)\frac{x_{1}^{\alpha (s_1+\cdots+s_n)}}{\prod_{i=1}^{n}s_{i}}\frac{(a_{i}+r)!(\frac{1+a_{i}+r}{\alpha}\cdots\frac{\alpha+a_{i}+r}{\alpha})_{s_1+\cdots+s_n}}{(a_{i})!(\frac{1+a_{i}}{\alpha}\cdots\frac{\alpha+a_{i}}{\alpha})_{s_1+\cdots+s_n}}\\
&=x_{1}^{a_{i}} (1+a_{i})_{r}\,F_{{l+\alpha}:\,m_1;\dots;\,m_{n}}^{p+\alpha:\,q_1;\dots;\,q_{n}}\Big[^{\frac{1+a_{i}+r}{\alpha},\dots, \frac{\alpha+a_{i}+r}{\alpha}, (a_p)\,\,:\,\,(b^{(1)}_{q_1})\,;\dots;\,(b^{(n)}_{q_{n}});}_{\,\frac{1+a_{i}}{\alpha},\dots, \frac{\alpha+a_{i}}{\alpha},\, (\alpha_l) \,:\,(\beta^{(1)}_{m_1})\,;\dots;\,(\beta^{(n)}_{m_{n}});} \,x_1^{\alpha},\dots, x_1^{\alpha}\Big].
\end{align*}
Applying the generalized Leibnitz formula for differentiation of the product of the following two functions, gives
\begin{align*}
&D_{x_{1}}^{r}\{x_{1}^{a_{i}+r}F_{{l}:\,m_1;\dots;\,m_{n}}^{p:\,q_1;\dots;\,q_{n}}\Big[^{(a_p)\,\,:\,\,(b^{(1)}_{q_1})\,;\dots;\,(b^{(n)}_{q_{n}});}_{(\alpha_l) \,:\,(\beta^{(1)}_{m_1})\,;\dots;\,(\beta^{(n)}_{m_{n}});} \,x_1^{\alpha},\dots, x_1^{\alpha}\Big]\}\\
&=\sum_{k=0}^{r}{r\choose k}\,D_{x_{1}}^{r-k}\{x_{1}^{a_{i}}\}\,D_{x_{1}}^{k}\{x_{1}^{r}F_{{l}:\,m_1;\dots;\,m_{n}}^{p:\,q_1;\dots;\,q_{n}}\Big[^{(a_p)\,\,:\,\,(b^{(1)}_{q_1})\,;\dots;\,(b^{(n)}_{q_{n}});}_{(\alpha_l) \,:\,(\beta^{(1)}_{m_1})\,;\dots;\,(\beta^{(n)}_{m_{n}});} \,x_1^{\alpha},\dots, x_1^{\alpha}\Big]\}\\
&=\sum_{k=0}^{r}{r\choose k} \frac{(-1)^{r+k}(-a_{i})_{r}(-r)_{k}}{(1+a_{i}-r)_{k}}\,x_{1}^{a_{i}}\\
&\qquad\times \,F_{{l+\alpha}:\,m_1;\dots;\,m_{n}}^{p+\alpha:\,q_1;\dots;\,q_{n}}\Big[^{\,\,\,\,\,\frac{1+r}{\alpha},\dots, \frac{\alpha+r}{\alpha},\,\, (a_p)\,\,:\,\,(b^{(1)}_{q_1})\,;\dots;\,(b^{(n)}_{q_{n}});}_{\frac{1+r-k}{\alpha},\dots, \frac{\alpha+r-k}{\alpha}, (\alpha_l) \,:\,(\beta^{(1)}_{m_1})\,;\dots;\,(\beta^{(n)}_{m_{n}});} \,x_1^{\alpha},\dots, x_1^{\alpha}\Big].
\end{align*}
Equating the above two identities,  get (\ref{1s13}). The other identity (\ref{m2s14}) can be proved  analogously.
\section{Finite summation formulas of generalized  Kamp\'e de F\'eriet series by rearrangement}
\begin{theorem}
The following finite summation formulas of generalized  Kamp\'e de F\'eriet series hold true:
\begin{align}
&\sum_{k=0}^{r}\frac{(-1)^{k}(a_{i+1})_{k}}{(a_{i+1}-a_{i}-r+1)_{k}}\,\,F_{{l}:\,m_1;\dots;\,m_{n}}^{p:\,q_1;\dots;\,q_{n}}\Big[^{a_{i+1}+k,\,(a_{p}^{i+1})\,:\,\,(b^{(1)}_{q_1})\,;\dots;\,(b^{(n)}_{q_{n}});}_{ \,\,\,\,\,\,\,\,\,\,(\alpha_l) \,\,\,\,\,\,\,\,\,\,\,\,\,\,\,\,:\,(\beta^{(1)}_{m_1})\,;\dots;\,(\beta^{(n)}_{m_{n}});} \,\,{x_1},\dots,{x_n}\Big]\notag\\
&=\,\frac{(a_{i})_{r}}{(a_{i}-a_{i+1})_{r}}F_{{l}:\,m_1;\dots;\,m_{n}}^{p:\,q_1;\dots;\,q_{n}}\Big[^{a_i+r, (a_p^{i})\,\,:\,\,(b^{(1)}_{q_1})\,;\dots;\,(b^{(n)}_{q_{n}});}_{\,\,\,\,\,\,\,\,\,\,\,\, (\alpha_l)\,\,\,\,\,\,:\,(\beta^{(1)}_{m_1})\,;\dots;\,(\beta^{(n)}_{m_{n}});} \,{x_1},\dots,{x_n}\Big],\label{m2s15}
\end{align}
where $i=1,\dots, p-1;$
\begin{align}
&\sum_{k=0}^{r}\frac{(-1)^{k}(b^{(1)}_{i+1})_{k}}{(b^{(1)}_{i+1}-b^{(1)}_{i}-r+1)_{k}}\,\,F_{{l}:\,m_1;\dots;\,m_{n}}^{p:\,q_1;\dots;\,q_{n}}\Big[^{ (a_p)\,\,:b^{(1)}_{i+1}+k, (b^{(1), i}_{q_1})\,;\dots;\,(b^{(n)}_{q_{n}});}_{ (\alpha_l) \,: \,\,\,\, (\beta^{(1)}_{m_1})\,;\dots;\,(\beta^{(n)}_{m_{n}});} \,\, x_{1}, \dots, {x_n}\Big]\notag\\
&=\,\frac{(b^{(1)}_{i})_{r}}{(b^{(1)}_{i}-b^{(1)}_{i+1})_{r}}F_{{l}:\,m_1;\dots;\,m_{n}}^{p:\,q_1;\dots;\,q_{n}}\Big[^{ (a_p)\,\,: b_{i}^{(1)}+r , (b^{(1), i}_{q_1})\,;\dots;\,(b^{(n)}_{q_{n}});}_{ (\alpha_l) \,: (\beta^{(1)}_{m_1})\,;\dots;\,(\beta^{(n)}_{m_{n}});} \,{x_1},  \dots, {x_n}\Big],\label{m2s16}
\end{align}
where $i=1,\dots, q_1-1.$
\end{theorem}
{\bf Proof:} Put the definition of generalized  Kamp\'e de F\'eriet series, the left hand side of (\ref{m2s15}) can be expressed as
\begin{align*}
&\sum_{s_1,\dots, s_n=0}^{\infty}\wedge(s_1,\dots, s_n) \prod_{i=1}^{n}\frac{x_{i}^{s_{i} }}{s_{i}!}\,_{2}F_{1}(-r, a_{i+1}+r; a_{i+1}-a_{i}-r+1; 1)
\\
&=\frac{(a_{i})_{r}}{(a_{i}-a_{i+1})_{r}}\,F_{{l}:\,m_1;\dots;\,m_{n}}^{p:\,q_1;\dots;\,q_{n}}\Big[^{a_i+r, (a_p^{i})\,\,:\,\,(b^{(1)}_{q_1})\,;\dots;\,(b^{(n)}_{q_{n}});}_{\,\,\,\,\,\,\,\,\,\,\,\, (\alpha_l)\,\,\,\,\,\,:\,(\beta^{(1)}_{m_1})\,;\dots;\,(\beta^{(n)}_{m_{n}});} \,{x_1},\dots,{x_n}\Big],
\end{align*}
where, we have used the Vandermonde's therorem
\begin{align*}
_{2}F_{1}(-r, a; b; 1)=\frac{(b-a)_{r}}{(b)_{r}},
\end{align*}
in the above equality. The other  identity  (\ref{m2s16}) can be proved in an analogous manner. 

\section{Conclusion}
Author established several finite summation formulas involving the generalized  Kamp\'e de F\'eriet series and remark that by specializing the parameters in generalized  Kamp\'e de F\'eriet series, author can deduce summation formulas for the generalized  Lauricella  functions \cite{SK, SM}  as well as  confluent  forms of Lauricella series in $n$ variables $\Phi_{2}^{(n)}$,  $\Psi_{2}^{(n)}$, $\Phi_{D}^{(n)}$, $\Xi_{1}^{(n)}$ and  $\Phi_{3}^{(n)}$, \cite{SK}. 
For example, specializing the parameters in (\ref{n2s2})  author get the  finite summation formulas for  $F_{B}^{(n)}$ and $\Xi_{1}^{(n)}$:
\begin{align}
&\sum_{k=0}^{r}{r\choose k}\frac{(b_1)_{k}}{(c)_{k}}x_{1}^{k}\,F_{B}^{(n)}(a_1+k, b_{1}+k, b_2,\dots,  b_{n}; c+k; {x_1}, \dots, {x_n} )\notag\\
&= F_{B}^{(n)}(a_1+r, b_{1}, \dots, b_{n}; c; {x_1}, \dots, {x_1} );
\end{align}
\begin{align}
&\sum_{k=0}^{r}{r\choose k}\frac{(b_1)_{k}}{(c)_{k}}x_{1}^{k}\,\Xi_{1}^{(n)}(a_1+k, a_2,\dots, a_n,  b_{1}+k, b_2,\dots,  b_{n-1}; c+k; {x_1}, \dots, {x_n} )\notag\\
&= \Xi_{1}^{(n)}(a_1+r, b_{1}, \dots, b_{n-1}; c; {x_1}, \dots, {x_1} );
\end{align}
Again, specializing the parameters in (\ref{1s11}) author obtain the following summation formula for  $F_{D}^{(n)}$ and $\Phi_{D}^{(n)}$:
\begin{align}
&\sum_{k=0}^{r}\frac{(-r)_{k}}{(1-2r)_{k}}\,F_{D}^{(n)}(1+r, b_{1},\dots,  b_{n}; 1;\frac{1}{x_1}, \dots, \frac{1}{x_1} )\notag\\
&=2 \,F_{D}^{(n)}(1+2r, b_{1}, \dots, b_{n}; 1+r;\frac{1}{x_1}, \dots, \frac{1}{x_1} );
\end{align}
\begin{align}
&\sum_{k=0}^{r}\frac{(-r)_{k}}{(1-2r)_{k}}\,\Phi_{D}^{(n)}(1+r, b_{1},\dots,  b_{n-1}; 1;\frac{1}{x_1}, \dots, \frac{1}{x_1} )\notag\\
&=2 \,\Phi_{D}^{(n)}(1+2r, b_{1}, \dots, b_{n-1}; 1+r;\frac{1}{x_1}, \dots, \frac{1}{x_1} ).
\end{align}
The details can be worked out by the reader.

\end{document}